\documentclass[12pt, twoside, reqno]{article}
\usepackage{amsmath,amsthm}
\usepackage{amssymb}
\usepackage{enumitem}

\newtheorem{theorem}{Theorem}[section]
\newtheorem{corollary}[theorem]{Corollary}
\newtheorem{lemma}[theorem]{Lemma}

\theoremstyle{definition}

\newtheorem{conjecture}[theorem]{Conjecture}

\numberwithin{equation}{section}

\begin{document}

\baselineskip=17pt


\title{The number of primes not in a numerical semigroup}

\author{Yong-Gao Chen\footnote{Corresponding author} and Hui Zhu\\
\small School of Mathematical Sciences, Nanjing Normal University \\
\small  Nanjing  210023,  P.R. China\\
\small ygchen@njnu.edu.cn, zh200109150528@163.com}
\date{}

\maketitle

\begin{abstract} For two coprime positive integers $a$ and $b$,
let $\pi^* (a, b)$ be the number of primes that cannot be
represented as $au+bv$, where $u$ and $v$ are nonnegative
integers. It is clear that $\pi^* (a, b)\le \pi (ab-a-b)$, where
$\pi (x)$ denotes the number of primes not exceeding $x$. In this
paper, we prove that $\pi^* (a, b)\ge 0.04\pi (ab-a-b)$ and pose
following conjecture: $\pi^* (a, b)\ge \frac 12 \pi (ab-a-b)$.
This conjecture is confirmed for $1\le a\le 10$.
\end{abstract}

\renewcommand{\thefootnote}{}



\renewcommand{\thefootnote}{\arabic{footnote}}
\setcounter{footnote}{0}

\section{Introduction}

For two coprime positive integers $a$ and $b$, let $\langle
a,b\rangle$ be the additive semigroup generated by $a$ and $b$,
i.e. $\langle a, b\rangle$ is the set of all integers of the form
$au+bv$, where $u$ and $v$ are nonnegative integers. In 1882 and
1883, Sylvester \cite{Sylvester1882} \cite{Sylvester1883} studied
integers of the form $au+bv$. It is well known that every integer
$>ab-a-b$ is in $\langle a, b\rangle$.  In 2020, Ram\'irez
Alfons\'in and Ska\l ba \cite{Ramiez2020} considered those primes
which are in $\langle a, b\rangle$ and less than $ab-a-b$. They
\cite{Ramiez2020} gave a lower bound for the number of such primes
and conjectured that the number is about $\frac 12 \pi (ab-a-b)$
as $\min\{ a, b\} \to +\infty$, where $\pi (x)$ denotes the number
of primes not exceeding $x$. In 2023, Ding \cite{Ding2023} showed
that the conjecture is true for almost all $a, b$ with $\gcd (a,
b)=1$.  Ding, Zhai and Zhao \cite{DingZhaiZhaoarXiv2023} confirmed
this conjecture.

In this paper, we consider those primes that are not in $\langle
a, b\rangle$. Let $S(a,b)=ab-a-b$ and  let $\pi^* (a, b)$ be the
number of such primes. It is clear that $\pi^* (a, b)\le \pi
(ab-a-b)$, $\pi^* (1, b)=0$, $\pi^* (2, 3)=0$ and $\pi^* (2,
b)=\pi (b-2)-1$ for $b>3$ and $2\nmid b$. In this paper, the
following results are proved.

\begin{theorem}\label{thm2} For any integer $a\ge 3$, there exists $C_a$ such that if
 $b$ is an integer with $b\ge C_a$ and $\gcd (a,b)=1$, then
\begin{equation}\label{e5}\pi^*
(a,b)>\left( \frac 12 +\frac 1{2(a-1)}\right) \frac {S(a,b)}{\log
S(a,b)} .\end{equation} Furthermore, the factor $\frac 12 +\frac
1{2(a-1)}$ is the best possible.

In particular, if $3\le a\le 1200$, then $C_a =50a^2$ is
admissible.
\end{theorem}

\begin{theorem}\label{thm1} For $b>a\ge 1$ with $\gcd (a, b)=1$, we
have
\begin{equation*}\pi^*
(a,b)\ge 0.04 \pi (S(a,b)) .\end{equation*}
\end{theorem}

From the proof of Theorem  \ref{thm1}, one may see that $0.04$ can be improved.  We pose two conjectures here.

\begin{conjecture} \label{coj2} For $b>a\ge 3$ with $\gcd (a, b)=1$ and $(a, b)\not= (3, 4)$, $(3, 5)$, $(3,7)$,  we
have
$$\pi^* (a,b)> \left( \frac 12 +\frac 1{2(a-1)}\right) \frac {S(a,b)}{\log S(a,b)}.$$ \end{conjecture}

\begin{conjecture}\label{coj1} For $b>a\ge 1$ with $\gcd (a, b)=1$, we
have
$$\pi^* (a,b)\ge \frac 12 \pi (S(a,b))$$
and the equality holds if and only if one of the following four
cases holds: (1) $a=1$; (2) $a=2$ and $b=3$; (3) $a=2$ and $b=5$;
(4) $a=3$ and $b=5$.
 \end{conjecture}

\begin{theorem}\label{thm3} Conjectures \ref{coj2} and \ref{coj1} are both true for $a\le 10$.  \end{theorem}

\section{Proofs }

In the following, let $\mathcal{P}$ be the set of all primes and
$S=S(a,b)$ for simplicity. Let $\pi (x; m, l)$ denote the number
of primes $p\le x$ with $p\equiv l\pmod{m}$ and let $\varphi (n)$
be the Euler totient function. For any positive integer $a$, let
$\omega (a)$ denote the number of distinct prime factors of $a$.
We begin with the following lemmas.

\begin{lemma}\cite[Theorem 2]{MontgomeryVaughan1973} \label{lem7}
Let $k,l$ be two integers and let $x,y$ be two positive real
numbers with $1\le k<y$. Then
$$\pi (x+y; k,l)-\pi (x; k,l)<\frac{2y}{\varphi (k) \log (y/k)}.$$
\end{lemma}

\begin{lemma} (Siegel--Walfisz Theorem, see \cite[Chapter 22]{Davenport1980}) \label{pi(x;m,l)=}
Let $A$ be a positive real number and $m,l$ two coprime
integers with $1\le m\le (\log x)^A$. Then there is a positive
constant $D$ depending only on $A$ such that
$$\pi (x;m,l)=\frac 1{\varphi (m)} \int_2^x \frac 1{\log t} dt +O(x \exp (-D \sqrt{\log x})), $$
uniformly in $m$.
\end{lemma}

\begin{corollary} \label{cor1} Let $m, l$ be two coprime integers. Then there is $X(m)$ such that
for all $x\ge X(m)$,
$$\frac{x}{\varphi (m) \log x}< \pi (x;m,l)< \frac{x}{\varphi (m) \log x} \left( 1+\frac 5{2\log x}\right).$$
\end{corollary}

\begin{proof}
By Lemma \ref{pi(x;m,l)=} with $A=1$,
\begin{align*}\pi (x;m,l)&=\frac 1{\varphi (m)} \int_2^x \frac 1{\log t} dt +O(x \exp (-D \sqrt{\log
x}))\\
&=\frac 1{\varphi (m)} \frac x{\log x}+\frac 1{\varphi (m)} \frac
x{(\log x)^2}\\
&\ \ +\frac 1{\varphi (m)} \int_2^x \frac 2{(\log t)^3} dt +O(x
\exp (-D \sqrt{\log
x}))\\
&>\frac 1{\varphi (m)} \frac x{\log x}+\frac 1{\varphi (m)} \frac
x{(\log x)^2}+O(x \exp (-D \sqrt{\log
x}))\\
&>\frac 1{\varphi (m)} \frac x{\log x}
\end{align*}
and
\begin{align*}\pi (x;m,l)&=\frac 1{\varphi (m)} \int_2^x \frac 1{\log t} dt +O(x \exp (-D \sqrt{\log
x}))\\
&=\frac 1{\varphi (m)} \frac x{\log x}+\frac 1{\varphi (m)} \frac
x{(\log x)^2}\\
&\ \ +\frac 1{\varphi (m)} \int_2^x \frac 2{(\log t)^3} dt +O(x
\exp (-D \sqrt{\log
x}))\\
&<\frac 1{\varphi (m)} \frac x{\log x}+\frac 5{2\varphi (m)} \frac
x{(\log x)^2}
\end{align*}
for all $x\ge X(m)$, where $X(m)$ is a constant depending only on
$m$.

This completes the proof of Corollary \ref{cor1}.
\end{proof}

\begin{lemma}\cite[Corollary 1.6]{Illinois2018} \label{pi(x;m,l)estimatesmall}
Let $m, l$ be two coprime integers with $1\le m\le 1200$. Then for all $x\ge 50m^2$,
$$\frac{x}{\varphi (m) \log x}< \pi (x;m,l)< \frac{x}{\varphi (m) \log x} \left( 1+\frac 5{2\log x}\right).$$
\end{lemma}

\begin{lemma}\cite{RosserSchoenfeld1962} \label{pi(x)estimate} We
have
$$\pi (x)>\frac x{\log x},\quad x\ge 17$$
and $$\pi (x)<\frac x{\log x} \left( 1+\frac{3}{2\log
x}\right),\quad x>1.$$
\end{lemma}

\begin{proof}[Proof of Theorem \ref{thm2}] For $3\le a\le 1200$,
let $C_a=50a^2$. For $a>1200$, let $C_a=X(a)+1$, where $X(a)$ is
as in Corollary \ref{cor1}. Let $b$ be an integer with $b\ge C_a$
and $\gcd (a,b)=1$. By $\gcd (a,b)=1$,  we know that $b\ge 50a^2$
implies $b-1\ge 50a^2$. Thus, $b-1\ge 50a^2$ for $3\le a\le 1200$
and $b-1\ge X(a)$ for $a>1200$.

In the following, $u,v$ always denote integers. It is easy to see
that an integer $n$ cannot be represented as $n=au+bv$ $(u\ge 0,
v\ge 0)$ if and only if $n$ can be represented as $n=au+bv$ $(u<0,
0\le v\le a-1)$. Hence
\begin{align*}\pi^* (a, b)&=| \{ p\in \mathcal{P} : p=au+bv, u<0,
0\le v\le a-1\} |\\
&=\sum_{v=1}^{a-1} |\{ p\in \mathcal{P} : p<bv, p\equiv bv\pmod{a}
\} |\\
&\ge \pi (b-1; a, b)+\sum_{\substack{v=2\\ \gcd (v,a)=1 } }^{a-2}
\pi (bv; a, bv) +\pi (S; a, b(a-1)).
\end{align*}
 It follows from Lemma
\ref{pi(x;m,l)estimatesmall} and Corollary \ref{cor1} that
\begin{align*}\pi^* (a, b)&\ge  \frac 1{\varphi (a)} \frac{b-1}{\log (b-1)}+ \sum_{\substack{v=2\\ \gcd
(v,a)=1 }}^{a-2}  \frac 1{\varphi (a)} \frac{bv}{\log (bv)}+ \frac
1{\varphi (a)} \frac{S}{\log S}\\
&\ge \frac 1{\varphi (a)} \frac{b-1}{\log S}
+\sum_{\substack{v=2\\ \gcd (v,a)=1 } }^{a-2}  \frac 1{\varphi
(a)} \frac{bv}{\log S}+ \frac
1{\varphi (a)} \frac{b(a-1)-a}{\log S}\\
&=\frac b{\varphi (a)\log S} \sum_{\substack{v=1\\ \gcd (v,a)=1}
}^{a-1} v -\frac {a+1}{\varphi (a)\log S}.
\end{align*}
In the above arguments, for $a=3$ we appoint
$$\sum_{\substack{v=2\\ \gcd (v,a)=1} }^{a-2} \pi
(bv; a, bv)=\sum_{\substack{v=2\\ \gcd (v,a)=1 } }^{a-2}  \frac
1{\varphi (a)} \frac{bv}{\log S}=0.$$ Since
$$\sum_{\substack{v=1\\ \gcd (v,a)=1} }^{a-1} v=\sum_{\substack{v=1\\ \gcd (v,a)=1} }^{a-1}
(a-v),$$ it follows that
$$\sum_{\substack{v=1\\ \gcd (v,a)=1} }^{a-1} v=\frac 12 a\varphi (a).$$
Hence
\begin{align*}\pi^* (a, b)&\ge
\frac {ab}{2\log S} -\frac
{a+1}{\varphi (a)\log S}\\
&=\frac{a}{2(a-1)}\frac{b(a-1)}{\log S} -\frac
{a+1}{\varphi (a)\log S}\\
&=\frac{a}{2(a-1)}\frac{S+a}{\log S} -\frac {a+1}{\varphi (a)\log
S}\\
&>\left( \frac 12 +\frac 1{2(a-1)}\right)  \frac S{\log S}.
\end{align*}
In the last inequality, we use the fact that $\varphi (a) a^2\ge
2a^2>2(a^2-1)$.

Similarly, we have
\begin{align*}\pi^* (a, b)
&=\sum_{v=1}^{a-1} |\{ p\in \mathcal{P} : p<bv, p\equiv bv\pmod{a}
\} |\\
&\le \sum_{\substack{v=1\\ \gcd (v,a)=1 } }^{a-1} \pi (bv; a, bv)
+\sum_{\substack{v=1\\ \gcd (v,a)>1 }}^{a-1} |\{ p\in \mathcal{P}
: p<bv, p\equiv bv\pmod{a} \} |.
\end{align*}
It is clear that
\begin{align*}&\sum_{\substack{v=1\\ \gcd (v,a)>1 }}^{a-1} |\{ p\in
\mathcal{P} : p<bv, p\equiv bv\pmod{a} \} |\\
=&\sum_{v=1}^{a-1} |\{ p\in \mathcal{P}
: p=\gcd (v, a), p\equiv bv\pmod{a} \} |\\
=&\sum_{\substack{p\mid a\\ p\in \mathcal{P}}}
\sum_{\substack{v=1\\ bv\equiv p\pmod{a}\\ \gcd (v, a)=p}}^{a-1} 1
\le   \sum_{\substack{p\mid a\\ p\in \mathcal{P}}}
\sum_{\substack{v=0\\ bv\equiv p\pmod{a}}}^{a-1}
1\\
= & \sum_{\substack{p\mid a\\ p\in \mathcal{P}}} 1=\omega (a).
\end{align*}
Therefore,
\begin{align*}\pi^* (a, b)\le \sum_{\substack{v=1\\ \gcd (v,a)=1 } }^{a-1} \pi
(bv; a, bv)+\omega (a).\end{align*} By using the upper bounds in
Lemma \ref{pi(x;m,l)estimatesmall} and Corollary \ref{cor1}, we
have
\begin{align*}\pi^* (a, b)&\le  \sum_{\substack{v=1\\ \gcd (v,a)=1 }}^{a-1} \frac 1{\varphi
(a)} \frac{bv}{\log (bv)}\left( 1+\frac 5{2\log (bv)}\right) +\omega (a)\\
&\le \frac b{\varphi (a) \log b} \left( 1+\frac 5{2\log
b}\right) \sum_{\substack{v=1\\ \gcd (v,a)=1 }}^{a-1} v+\omega (a)\\
&= \frac 12 \frac {ab}{\log b} \left( 1+\frac 5{2\log
b}\right) +\omega (a)\\
&=\left( \frac 12 +\frac 1{2(a-1)}\right) \frac {(a-1)b}{\log b}
\left( 1+\frac 5{2\log b}\right) +\omega (a).\end{align*} Given
$a$. When $b\to +\infty$, we have
$$\frac {(a-1)b}S\to 1,\quad \frac {\log S}{\log
b} \to 1, \quad 1+\frac 5{2\log b} \to 1.$$ It follows that
\begin{align*}\pi^* (a, b)\le  \left( \frac 12 +\frac 1{2(a-1)} +o(1)\right) \frac{S}{\log S}
.\end{align*} Therefore, $\frac 12 +\frac 1{2(a-1)}$ in \eqref{e5}
is the best possible.

This completes the proof of Theorem \ref{thm2}.
\end{proof}

\begin{proof}[Proof of Theorem \ref{thm1}] Since $S(1, b)=-1$ and $S(2, b)=b-2$, it
follows that $\pi^* (1, b)\ge 0=\pi (S(1, b))$, $\pi^* (2, 3)=
0=\pi (S(2, 3))$ and for $b\ge 5$ with $2\nmid b$,  $$\pi^* (2,
b)=\pi (b-2)-1\ge \frac 12\pi (b-2)=\frac 12 \pi (S(2, b)).$$  In the following, we always
assume that $a\ge 3$. For $0< \delta \le 1$, we have
\begin{equation}\label{w1}\pi^* (a, b)\ge \pi (\delta S)-\# \{ p\in \mathcal{P} : 0< p\le \delta S, p\in
\langle a, b\rangle\} . \end{equation}
Let $p\in  \mathcal{P}$ such that $p= au+bv\le \delta S$ for some
$u,v\in \mathbb{Z}_{\ge 0}$. Then $\gcd (a, v)\mid p$ and $v<\delta a$ by $S=ab-a-b<ab$. If $v=0$, then $p=a$. If $v\ge 1$, then by $p=
au+bv\ge b>a\ge \gcd (a, v)$, we have $\gcd (a, v)=1$. Hence
\begin{align*}&\# \{ p\in \mathcal{P} : 0< p\le \delta S, p\in \langle a, b\rangle  \} \\
&\le \sum_{\substack{1\le v\le \delta a\\ \gcd (a,v)=1}} \# \{ p\in
\mathcal{P} : 0< p\le \delta S, p\equiv
bv\pmod{a} \} +1\nonumber\\
&= \sum_{\substack{1\le v\le \delta a\\ \gcd (a,v)=1}} \pi (\delta S; a,
bv) +1.
\end{align*}

If  $a<\delta S$, then by Lemma \ref{lem7},
\begin{align*}
 \pi (\delta S; a, bv) <  \frac{2\delta S}{\varphi
(a) \log (\delta S/a)} = \frac {2\delta }{\varphi (a)}\left( 1-\frac{\log a}{\log \delta S}\right)^{-1} \cdot \frac{S}{\log \delta S}.
\end{align*}
Hence
\begin{align*}&\# \{ p\in \mathcal{P} : 0< p\le \delta S, p\in \langle a, b\rangle  \} \nonumber\\
& < \frac {2\delta }{\varphi (a)}\Big( \sum_{\substack{1\le v\le \delta a\\
\gcd (a,v)=1}} 1\Big) \Big( 1-\frac{\log a}{\log
\delta S}\Big)^{-1}\cdot \frac{S}{\log \delta S} +1.\end{align*}
If $\delta S\ge 17$, then by Lemma \ref{pi(x)estimate},
$$\pi (\delta S)>\frac{\delta S}{\log \delta S}.$$
It follows from \eqref{w1} that if $\delta S\ge 17$ and  $a<\delta S$, then
\begin{align*}\pi^* (a, b)&> \frac{\delta S}{\log \delta S}-\frac {2\delta }{\varphi (a)}\Big( \sum_{\substack{1\le v\le \delta a\\
\gcd (a,v)=1}} 1\Big) \Big( 1-\frac{\log a}{\log
\delta S}\Big)^{-1}\cdot \frac{S}{\log \delta S} -1\\
&=\left( 1-\frac {2 }{\varphi (a)}\Big( \sum_{\substack{1\le v\le \delta a\\
\gcd (a,v)=1}} 1\Big) \Big( 1-\frac{\log a}{\log
\delta S}\Big)^{-1}- \frac{\log \delta S}{\delta S} \right) \frac{\delta S}{\log \delta S}.\end{align*}
If $\delta <e^{-3/2}=0.2231\cdots $, then by Lemma \ref{pi(x)estimate},
\begin{align*}
\pi (S)&<\left( 1+\frac{3}{2\log S}\right) \frac{S}{\log S}\\
&<\left( 1-\frac{\log \delta }{\log S}\right) \frac{S}{\log S}\\
&<\left( 1+\frac{\log \delta }{\log S}\right)^{-1} \frac{S}{\log S}\\
&=\frac{S}{\log \delta S}.
\end{align*}
It follows from \eqref{w1} that if $\delta S\ge 17$ and  $a<\delta S$, then
\begin{align}\label{w2}\pi^* (a, b)>\Delta (\delta , a, S) \pi (S),
\end{align}
where
\begin{align*}\Delta (\delta , a, S) =\left( 1-\frac {2 }{\varphi (a)}\Big( \sum_{\substack{1\le v\le \delta a\\
\gcd (a,v)=1}} 1\Big) \Big( 1-\frac{\log a}{\log
\delta S}\Big)^{-1}- \frac{\log \delta S}{\delta S} \right) \delta .
\end{align*}

Now we divide into the following cases:

{\bf Case 1:} $a>6\cdot 10^4$. Since
\begin{align*}\sum_{\substack{1\le v\le 0.1a\\
\gcd (a,v)=1}} 1&=\sum_{1\le v\le 0.1a} \sum_{d\mid a, d\mid v}
\mu (d)=\sum_{d\mid a} \mu (d)\sum_{1\le v\le 0.1a, d\mid v} 1=\sum_{d\mid a} \mu (d)\left\lfloor \frac{0.1a}{d}\right\rfloor \\
&\le \sum_{d\mid a} \mu (d)\frac{0.1a}{d} + \sum_{d\mid a}|\mu
(d)|=0.1 \varphi (a) +2^{\omega (a)},
\end{align*}
it follows that
$$\frac {1}{\varphi (a)}\sum_{\substack{1\le v\le 0.1a\\
\gcd (a,v)=1}} 1\le 0.1+\frac{2^{\omega (a)}}{\varphi (a)}.$$ For
a prime $p\ge 2$ and an integer $k\ge 1$, we have
$$\frac{2p^{k/2}}{p^{k-1}(p-1)} =\frac{2p^{1-k/2}}{p-1}\le
\frac{2\sqrt{p}}{p-1}.$$ It is easy to see that
$$\frac{2\sqrt{p}}{p-1}>1, \quad p\in \mathcal{P}$$ if and only if $p\in \{ 2, 3, 5 \} $.
Hence
$$\frac{2^{\omega (a)}}{\varphi (a)}=\prod_{p^k\| a}
\frac{2}{p^{k-1}(p-1)}\le \frac{2\sqrt{2}}{2-1}\cdot
\frac{2\sqrt{3}}{3-1}\cdot \frac{2\sqrt{5}}{5-1} \prod_{p^k\| a}
p^{-k/2}= \frac{\sqrt{30}}{\sqrt a}.$$ where the products are
taken over all prime powers $p^k$ such that $p^k\mid a$ and
$p^{k+1}\nmid a$. It follows from $a>6\cdot 10^4$ that
\begin{align*}\frac {1}{\varphi (a)}\sum_{\substack{1\le v\le 0.1a\\
\gcd (a,v)=1}} 1\le 0.1+\frac{\sqrt{30}}{\sqrt a}<0.122362.
\end{align*}
By $a>6\cdot 10^4$ and $ S=(a-1)(b-1)-1\ge (a-1)a-1$, we have $0.1 S\ge 17$, $a<0.1S$ and
\begin{align}\label{w3}\frac{\log a}{\log 0.1 S}\le \frac{\log a}{\log 0.1 (a^2-a-1)} <0.558438, \end{align} \begin{align}\label{w4}\frac{\log
0.1 S}{0.1 S}\le  \frac{\log 0.1 (a^2-a-1)}{0.1
(a^2-a-1)}<5.5\cdot 10^{-8} .\end{align} It follows that
\begin{align*}\Delta (0.1 , a, S)&>(1-2\times 0.122362 (1-0.558438)^{-1}-5.5\cdot 10^{-8})\times 0.1\\
&>0.0445 .\end{align*} By \eqref{w2}, $\pi^* (a, b)>0.0445 \pi
(S).$

{\bf Case 2:} $181\le a\le 6\cdot 10^4$. Let $\delta =0.0904$. By $ S=(a-1)(b-1)-1\ge (a-1)a-1$, we have $\delta S\ge 17$ and $a<\delta S$.
In view of \eqref{w2}, \eqref{w3}, \eqref{w4} and a calculation, we have
\begin{align*}\Delta (\delta , a, S)&>\Big( 1-\frac {2 }{\varphi (a)}\Big( \sum_{\substack{1\le v\le \delta a\\
\gcd (a,v)=1}} 1\Big) \Big( 1-\frac{\log a}{\log
(\delta h(a))}\Big)^{-1}- \frac{\log (\delta h(a))}{\delta h(a)} \Big) \delta\\
&>0.0401,\end{align*} where $h(a)=a^2-a-1$. By \eqref{w2}, $\pi^*
(a, b)>0.0401 \pi (S).$

{\bf Case 3:} $16\le a\le 180$. If $b\le 1000$ with $\gcd (a, b)=1$, then
a calculation shows that
\begin{align*}\pi^* (a, b)&\ge \pi (0.05 S)-\# \{ (u, v)\in \mathbb{Z}^2 :  u\ge 0,
v\ge 0, au+bv\le 0.05 S, au+bv\in \mathcal{P} \}\\
&> 0.0663 \pi (S).
\end{align*}
If $b>1000$, let $\delta =0.095$,  then $S=ab-a-b\ge 1000a-1001$,
$\delta S\ge 17$ and $a<\delta S$.  By a calculation,  we have
\begin{align*}\Delta (\delta , a, S)&>\Big( 1-\frac {2 }{\varphi (a)}\Big( \sum_{\substack{1\le v\le \delta a\\
\gcd (a,v)=1}} 1\Big) \Big( 1-\frac{\log a}{\log
(\delta g(a))}\Big)^{-1}- \frac{\log (\delta g(a))}{\delta g(a)} \Big) \delta\\
&>0.0425,\end{align*} where $g(a)=1000a-1001$. By \eqref{w2},
$\pi^* (a, b)>0.0425\pi (S).$

{\bf Case 4:} $3\le a\le 15$.
If $b\ge 181$,  then by Lemma \ref{pi(x)estimate} and $S=ab-a-b\ge 180a-181$, we have
\begin{align*}\pi^* (a, b)&\ge \pi (b-1)-1>\frac {b-1}{\log (b-1)}-1=\frac 1{a}
\frac{a(b-1)}{\log (b-1)}-1\\
&>\frac 1{a}\frac S{\log S}-1= \left(\frac 1{a}-\frac{\log
S}S\right) \frac S{\log S}\\
&>\left(\frac 1{a}-\frac{\log S}S\right) \left( 1+\frac 3{2\log
S}\right)^{-1} \pi (S)\\
&>\left(\frac 1{a}-\frac{\log (180a-181)}{180a-181}\right) \left(
1+\frac 3{2\log
(180a-181)}\right)^{-1} \pi (S)\\
&>0.05334 \pi (S).
\end{align*}
If $b\le 180$ with $\gcd (a, b)=1$, then a calculation shows that
\begin{align*}\pi^* (a, b)&=\pi (S)-\# \{ (u, v)\in \mathbb{Z}^2 :  u\ge 0,
v\ge 0, au+bv\le S, au+bv\in \mathcal{P} \}\\
&\ge  \frac 12 \pi (S).\end{align*}

This completes the proof of Theorem \ref{thm1}.
\end{proof}

\begin{proof}[Proof of Theorem \ref{thm3}]
For $a=1$, we have $S(1, b)=-1$ and $$\pi^* (1,b)=0= \frac 12 \pi
(S(1,b)).$$ For $a=2$ and $b=3,5$, we have
$$\pi^* (2,3)=0=\frac 12 \pi (S(2,3))$$
and
$$\pi^* (2,5)=1=\frac 12 \pi (S(2,5)).$$
For $a=2$ and $b\ge 7$ with $2\nmid b$, we have $S(2,b)=b-2\ge 5$.
It follows that $\pi (S(2,b))=\pi (b-2)\ge 3$ and $$\pi^*
(2,b)=\pi (b-2)-1\ge \frac 23 \pi (b-2)>\frac 12 \pi (S(2,b)).$$
Hence, Conjecture \ref{coj1} is true for $a=1,2$.

By Theorem \ref{thm2}, for $3\le a\le 10$ and $b\ge 50 a^2$ with
$\gcd (a, b)=1$ we have
\begin{equation}\label{eq1}\pi^* (a,b)>\left( \frac 12 +\frac 1{2(a-1)}\right) \frac
{S}{\log S}.\end{equation} For $3\le a\le 10$ and $a<b<50 a^2$
with $\gcd (a,b)=1$, a calculation shows that \eqref{eq1} holds
except for $(a, b)=(3, 4), (3, 5), (3, 7)$. In view of Lemma
\ref{pi(x)estimate},
\begin{equation}\label{eq2}\pi (S)<\frac {S}{\log S} \left( 1+\frac 3{2\log
S}\right) .\end{equation} If $S>e^{3(a-1)/2}$, then by
\eqref{eq2},
\begin{equation}\label{eq3a}\pi (S)<\frac {S}{\log S} \left( 1+\frac 1{a-1}\right)
.\end{equation} By \eqref{eq1} and \eqref{eq3a},
\begin{equation}\label{eq3}\pi^* (a,b)>\frac 12\pi (S).\end{equation}
If $b\ge \frac 1{a-1} e^{3(a-1)/2}+2$, then $b>7$ and
$$S(a,b)=(a-1)b-a\ge e^{3(a-1)/2}+2(a-1)-a>e^{3(a-1)/2}.$$ It
follows  that \eqref{eq3} holds. For $3\le a\le 8$ and  $b< \frac
1{a-1} e^{3(a-1)/2}+2$ with $\gcd (a,b)=1$, a calculation gives
\begin{align}\label{eq4}\pi^* (a,b)&= \pi (S) -\sum_{u=0}^{b-2}\sum_{v=0}^{a-2}|\{ (u,v) : au+bv\le S, au+bv\in \mathcal{P} \} |\nonumber\\
   &> \frac 12\pi (S)\end{align}
except for $(a, b)=(3,5)$. It is easy to see that
$$\pi^* (3,5) =2=\frac 12 \pi (S(3,5)).$$

For $a=9$, if $18595<S\le e^{3(a-1)/2}$, then \eqref{eq3a} holds
by a calculation and so \eqref{eq3} holds. If $S\le 18595$, then
$b\le 2325$. A calculation shows that \eqref{eq4} holds for
$a=9<b\le 2325$ with $\gcd (9,b)=1$.

For $a=10$, if $60180<S\le e^{3(a-1)/2}$, then \eqref{eq3a} holds
by a calculation and so \eqref{eq3} holds. If $S\le 60180$, then
$b\le 6687$. A calculation shows that \eqref{eq4} holds for
$a=10<b\le 6687$ with $\gcd (10,b)=1$.

This completes the proof of Theorem \ref{thm3}.
\end{proof}

\section*{Acknowledgments}

This work was supported by the National Natural Science Foundation
of China, Grant No. 12171243.

\end{document}